\input amstex
\documentstyle{amsppt}
\nopagenumbers \noindent {\bf Erratum for:
Inverse Spectral Theory for a Singular
Sturm-Liouville Operator on [0,1]} \vskip.1in
 \centerline {J.-C. Guillot and J. V. Ralston}

\vskip.2in We are grateful to Gerald Teschl for
calling our attention to an error in our paper
[1]. Our claim (p. 357, line 6) that the
solutions $\phi(x,\lambda)$ and $\tilde
\psi(x,\lambda)$ to the equation
$$(H-\lambda)u=_{\hbox{def.}}-{d^2u\over dx^2}+{2\over x^2}u+q(x)u-\lambda u=0\eqno{(1)}$$
are linearly independent for all $\lambda\in \Bbb
C$ is quite false (see [2, Remark 2.7]). In [3]
Kostenko, Sakhnovich and Teschl have given a
construction of a new solution which is entire in
$\lambda$ and linearly independent from $\phi$
for all $\lambda$.

For the purposes of [1] (and related subsequent
papers) one can construct $\tilde \psi$ by
replacing the inhomogeneous term $v$ in the
Volterra equation (1.10) by the function
$$w(x,\lambda)=\sqrt \lambda w_1^+(x\sqrt \lambda )=
3v-{i\over 3}\lambda^{3/2}u,$$ where $u$ and $v$
are the functions defined in formulas (1.4) and
(1.5) in [1].  Here $w_1^+$ is a Ricatti-Hankel
function (see [4, p.38]). With this modification
the extension of $\tilde \psi$ to $x>1$ as a
solution (1) becomes $w$. If for
$\lambda_0\in\Bbb R_+$ the functions
$\phi(x,\lambda_0)$ and $\tilde
\psi(x,\lambda_0)$ were linearly dependent,
$\tilde \psi(x,\lambda_0)$ would be a solution of
(1) for $0<x<\infty$ which vanishes
 at $x=0$ and is asymptotic to
$i\sqrt{\lambda_0}\exp(ix\sqrt{\lambda_0})$ as
$x\to\infty$. If the potential $q$ is
sufficiently well behaved that $\phi^\prime$
absolutely continuous on $[0,1]$ (for instance
$\int_0^1x|q(x)|dx<\infty$), then one has
$$0=\int_0^R[\overline{\tilde \psi}(H-\lambda_0)\tilde
\psi-\tilde \psi(H-\lambda_0)\overline{\tilde
\psi}]dx=[-\tilde \psi^\prime\overline{ \tilde
\psi}+
\overline{\tilde\psi^\prime}\tilde\psi](R)\to 2i
\lambda_0$$ as $R\to \infty$. Hence $\phi$ and
$\tilde \psi$ are linearly independent for
$\lambda$ on the positive real axis. This method
of constructing an \lq\lq irregular" solution to
(1) is given in [4, p. 373]. This reference was
cited in [1], but not followed closely enough.

The functions $\phi$ and $\tilde \psi$ are still
not necessarily linearly independent for all
$\lambda \in \Bbb C$, but that is not needed in
[1]. The first coordinate in the mapping $\Phi$
that introduces coordinates on the space of
potentials is $\int_0^1q(x)dx$, but the rest of
the coordinates are unchanged when one replaces
$q$ by $q+c,\ c\in \Bbb R$. Linear independence
of $\phi$ and $\tilde \psi$ at the eigenvalues
$\mu_n(q)$ is all that is used in [1], and in the
analysis of all but the first coordinate we can
assume that all of these eigenvalues are
positive.

The corrections to the remainder of [1] are as
follows. The estimates (1.12) and (1.13) hold for
the new $\tilde \psi$, but the exponential
factors in those estimates should be
$\exp(|\hbox{Im}\lambda^{1/2}|(2-x))$ instead of
$\exp(|\hbox{Im}\lambda^{1/2}|(1-x))$ -- this was
a typographic error in [1]. In (1.10), (1.13),
and (1.29) of [1] $v$ should be replaced by $w$.

In [5] the same changes are needed: one should
use constant multiples of the Ricatti-Hankel
functions $w_l^+$ as the inhomogeneous terms in
the integral equations for the second solutions
of

$$-{d^2u\over dx^2}+{l(l+1)\over x^2}u+q(x)u-\lambda u=0.$$
\newpage

\centerline{\bf References} \vskip.1in \noindent
[1] J.-C. Guillot and J. V. Ralston, Inverse
spectral theory for a singular Sturm-Liouville
operator on [0,1], J. Diff. Eq. {\bf 76}(1988),
353-373. \vskip.1in \noindent [2] A. Kostenko, A.
Sakhnovich and G. Teschl, Inverse eigenvalue
problems for perturbed spherical Schroedinger
operators, arXiv:1004.4175. \vskip.1in \noindent
[3] A. Kostenko, A. Sakhnovich and G. Teschl,
Weyl-Titchmarsh Theory for

\noindent Schroedinger operators with strongly
singular potentials, arXiv:1007.0136 \vskip.1in
\noindent [4] R. G. Newton, \lq\lq Scattering
Theory of Waves and Particles," McGraw-Hill, New
York, 1966. \vskip.1in \noindent [5] F. Serier,
The inverse spectral problem for radial
Schr\"odinger operators on [0,1], J. Diff. Eq.
{\bf 235}(2007), 101-126. \vskip.2in \noindent
(J.-C. Guillot) Centre de Math\'ematiques
Appliqu\'ees, UMR 7641, \'Ecole Polytech-

\noindent nique-CNRS, 91128 Palaiseau Cedex,
France.

\noindent {\it email:}
jean-claude.guillot\@polytechnique.edu \vskip.2in
\noindent (J.V. Ralston) Dept. of Math., UCLA,
Los Angeles, CA 90095, USA.

\noindent {\it email:} ralston\@math.ucla.edu
\end